\documentclass[11pt]{amsart}
\usepackage{amssymb}
\usepackage{amsmath}
\usepackage{enumerate}
\usepackage[all]{xypic}
\usepackage{hyperref}

\DeclareFontFamily{U}{rsfs}{\skewchar\font"7F}
\DeclareFontShape{U}{rsfs}{m}{n}{
	<-6> rsfs5
	<6-8> rsfs7
	<8-> rsfs10
	}{}
\DeclareMathAlphabet{\mathscr}{U}{rsfs}{m}{n}

\newtheorem{thm}{Theorem}[section]
\newtheorem{lem}[thm]{Lemma}
\newtheorem{cor}[thm]{Corollary}
\newtheorem{prop}[thm]{Proposition}
\newtheorem{conj}[thm]{Conjecture}

\theoremstyle{definition}

\newtheorem{rem}[thm]{Remark}
\newtheorem{defn}[thm]{Definition}

\newtheorem{ex}[thm]{Example}

\newtheorem{exs}[thm]{Examples}

\theoremstyle{remark}
\newtheorem{notation}[thm]{Notation}

\numberwithin{equation}{section}

\def\top{{\mathrm {top}}}
\def\red{{\mathrm {r}}}

\DeclareMathOperator{\Ind}{Ind}
\DeclareMathOperator{\Res}{Res}
\DeclareMathOperator{\Hom}{Hom}
\DeclareMathOperator{\Ext}{Ext}
\DeclareMathOperator*{\colim}{colim}
\DeclareMathOperator{\coker}{coker}

\def\KK{{K\!K}}

\def\Id{\mathrm{Id}}

\def\KKH{\KK^{\mathrm{naive}}}
\def\KKMN{\KK^{\top}}

\def\XG{{X \r G}}

\def\CI{\mathscr{C\!I}}
\def\LCI{\langle \CI \rangle}
\def\E{{\mathcal E}}
\def\pt{{{\{\mathrm {pt}}\}}}
\def\S{{\mathcal S}}

\def\F{{\mathcal F}}
\def\G{{\mathcal G}}

\def\Ab{{\mathbf {Ab}}}
\def\P{P}

\def\Co{{C_{0}}}

\def\x{{\otimes}}
\def\r{{\rtimes}}
\def\t{{\times}}

\def\ca{{$C^\ast$-algebra}}
\def\ga{{$G$-algebra}}
\def\gs{{$G$-space}}

\def\xga{{$X \r G$-algebra}}

\def\cont{{continuous}}

\def\BC{{\mathcal{BC}}}
\def\GD{{Going-Down }}
\def\lc{{locally compact}}

\def\Z{{\mathbb Z}}
\def\C{{\mathbb C}}

\mathchardef\ordinarycolon\mathcode`\: 
\def\vcentcolon{\mathrel{\mathop\ordinarycolon}} 
\providecommand*\coloneqq{\mathrel{\vcentcolon\mkern-1.2mu}=}

\begin{document}

\title[The Baum-Connes conjecture for $\KK$-theory]
{The Baum-Connes conjecture for $\KK$-theory}

\author[Otgonbayar Uuye]{Otgonbayar Uuye}
\address{Department of Mathematical Sciences,
University of Copenhagen, Denmark}
\email{otogo@math.ku.dk}

\begin{abstract}
We define and compare two bivariant generalizations of the topological $K$-group $K^\top(G)$ for a topological group $G$. We consider the Baum-Connes conjecture in this context and study its relation to the usual Baum-Connes conjecture. 
\end{abstract}

\maketitle

\section{Introduction}
$K$-theory has been one of the most successful tools for analyzing \ca s and $C^\ast$-dynamical systems. In this paper we consider the Baum-Connes conjecture, which proposes a way to compute the $K$-theory of a reduced crossed product algebra (see Section~\ref{Baum-Connes} for more details):
\begin{conj}[The Baum-Connes Conjecture with Coefficients]\label{BCC} Let $G$ be a \lc~ second-countable topological group. Then for any \ga~ $B$, the reduced assembly map $$\beta^B_\red: K^\top_{*}(G; B) \longrightarrow K_{*}(B \r_\red G)$$ is an isomorphism.
\end{conj}
If this is the case, we say that $G$ satisfies the Baum-Connes conjecture for $B$.

Counterexamples to the Baum-Connes conjecture were constructed by Higson, Lafforgue and Skandalis, building on ideas of Gromov, in \cite{HLS}. Nonetheless, the conjecture for $B = \C$ still stands and has profound applications to geometry and algebra.

In order to study the $\KK$-class of  $B \r_{\red} G$, we would like to generalize the conjecture to $\KK$-theory. This would allow, in particular, to determine the mod-$n$ $K$-theory of $B \r_{\red} G$.

The formulation of the Baum-Connes conjecture (with coefficients) given in \cite[Conjecture 9.6]{BCH} has a straightforward generalization to $\KK$-theory (cf.\ Conjecture~\ref{naive BC for KK}). However, one can easily see that while the right-hand-side of the conjecture is $\sigma$-additive in the first variable, the left-hand-side is {\em not}, in general. Hence this generalization of the conjecture to $\KK$-theory fails for ``trivial'' reasons. 

Meyer and Nest gave a reformulation of the Baum-Connes conjecture in \cite[Theorem 5.2]{MN-T}, using the notion of a Dirac morphism. Their approach yields another generalization of the conjecture to $\KK$-theory (cf.\ Conjecture~\ref{BCMN}), which behaves better in many respects. We remark that this generalization also has well-understood counter-examples (cf.\ Example \ref{strong BC examples}(\ref{strong BC for Sp})), but we believe it still serves as a useful tool in the study of the $\KK$-class of crossed product algebras.

In this paper, we compare the two approaches. In order to distinguish the two, we call the version based on \cite{BCH}, Conjecture \ref{naive BC for KK}, the {\em naive Baum-Connes conjecture for $\KK$-theory}, short for the naive generalization of the Baum-Connes conjecture to $\KK$-theory and the version based on \cite{MN-T}, Conjecture~\ref{BCMN}, simply, the {\em Baum-Connes conjecture for $\KK$-theory}. We often omit the ``for $\KK$-theory'' part.

Our main theorem is the following, see Theorem \ref{thm comparison} for the precise statement.
\begin{thm}\label{thm main intro} Let $B$ be a \ga. If the functor $\KK_{\ast}(A, -)$ commutes with colimits, then the two generalizations of the Baum-Connes conjecture to $\KK$-theory are equivalent for $(A, B)$.
\end{thm}

If $A$ satisfies the Universal Coefficient Theorem (cf.\ Theorem~\ref{UCT}) and has finitely generated $K$-theory, then $A$ satisfies the condition of Theorem~\ref{thm main intro} (cf.\ \cite[Theorem 7.13]{RS}). A particular example is the dimension-drop algebra $\mathbb I_{n}$, $n \ge 2$, of (\ref{dimension-drop}). Since the mod-$n$ $K$-theory of an algebra $D$ can be computed as 
	\begin{equation}K_{\ast}(D; \Z/n\Z) \cong \KK_{\ast}(\mathbb I_{n}, D),
	\end{equation}
(see \cite{DL}), we can consider the (naive) Baum-Connes conjecture for $(\mathbb I_{n}, B)$ as a Baum-Connes conjecture for $B$ in mod-$n$ $K$-theory. It follows from Theorem~\ref{thm main intro}, the two versions are equivalent. Moreover, they follow from the usual Baum-Connes conjecture:
	
\begin{thm}[Corollary \ref{UCT and BCMN} and Corollary \ref{UCT and naive BC}]\label{mod n} Let $B$ be a \ga~ for which $G$ satisfies the Baum-Connes conjecture (Conjecture \ref{BCC}). Then for any $A$ satisfying the UCT, $G$ satisfies the Baum-Connes conjecture for $(A, B)$. If in addition, $A$ has finitely generated $K$-theory, then $G$ satisfies the {\em naive} Baum-Connes conjecture for $(A, B)$.
\end{thm}

This is an immediate corollary of the treatment of UCT given in Section~\ref{section UCT}.

{\bf Acknowledgements.} This paper grew out of my master's thesis written at the University of Tokyo. I would like to thank my advisor Yasuyuki Kawahigashi for his generous support and constant encouragement, and J. Chabert, S. Echterhoff and N. Higson for many helpful discussions and comments. I am also very grateful to the referee, whose comments lead to substantial improvements in the presentation.

\section{Conventions}\label{sec Conventions}

Throughout the paper, we assume that topological groups and topological spaces are {\em second-countable, \lc~ and Hausdorff}, unless stated otherwise. Similarly, \ca s are tacitly assumed to be {\em separable}, with the obvious exceptions such as multiplier algebras. 
 
Let $G$ be a topological group and let $X$ be a topological space.  A {\em $G$-algebra} is a \ca~ equipped with a strongly continuous action of $G$. If $A$ is a \ga~ equipped with the trivial action of $G$, we often simply say that ``$A$ is a \ca". A {\em $\Co(X)$-algebra} is a \ca~equipped with a $\Co(X)$-action, that is, a {\em nondegenerate} $\ast$-homomorphism from $\Co(X)$ to the central multipliers of the algebra. Here $\Co(X)$ denote the \ca~ of continuous functions on $X$ vanishing at infinity.

Suppose that $X$ is a {\em $G$-space}, that is, $X$ is equipped with a continuous action of $G$. Then the algebra $\Co(X)$ is naturally a $G$-algebra via $(g \cdot f)(x) = f(g^{-1}x)$ for $g \in G$ and $f \in \Co(X)$. An {\em \xga~} is a $G$-$\Co(X)$-algebra such that the action of $\Co(X)$ is {\em $G$-equivariant}.

We say that $X$ is {\em $G$-compact} if the quotient $X/G$ is compact and {\em proper} if the map $X \t G \to X \t X$, $(x, g) \mapsto (x, gx)$ is proper. Note that \cite{BCH} considers a slightly different notion of properness; see \cite{BMP,Bi} for comparison. A \ga~is said to be {\em proper} if it can be obtained from an \xga~ with $X$ proper  by forgetting the $\Co(X)$-action. Note that for a proper algebra the reduced and full crossed products coincide. If $A$ is a \ga, following Kasparov, we write $A(X)$ for $\Co(X, A) = A \otimes \Co(X)$ and equip it with the diagonal action. 

Let $H$ be a closed subgroup of $G$. We denote the {\em restriction} (cf.\ \cite[Definition 3.1]{Ka}) and the {\em induction} (cf.\ \cite[Theorem 3.5]{Ka}) functors of Kasparov by $\Res^{G}_{H}: \KK^{X \r G} \to \KK^{X \r H}$ and $\Ind^{G}_{H}: \KK^{X \r H} \to \KK^{X \r G}$,  respectively.
\section{The Baum-Connes Conjecture for $\KK$-theory, attempt $1$}\label{Baum-Connes}
In this section, we consider the most simple-minded generalization of the Baum-Connes conjecture to $\KK$-theory and show why this is not the desired one.

\subsection{The naive Baum-Connes conjecture for $\KK$-theory}
Let $G$ be a topological group and let $\E{G}$ denote a {\em universal} proper space. (cf.\ \cite{BCH,BMP,KS2}). 

\begin{defn}\label{def top KK} Let $A$ be a \ca\footnote{Considered a $G$-algebra with the trivial action of $G$.} and let $B$ be a \ga. An inclusion $Y_{1} \subseteq Y_{2}$ of $G$-compact subsets of $\E{G}$  induces a natural map
	\begin{equation}
	\KK^{G}_{\ast}(A(Y_{1}), B) \to \KK^{G}_{\ast}(A(Y_{2}), B).
	\end{equation}
We define the \textit {naive topological $\KK$-groups of $(A, B)$} as 
	\begin{equation}
	\KKH_{\ast}(G; A, B) \coloneqq \colim_{\substack{Y \subseteq \E{G} \\\text{$G$-compact}}} \KK^G_{\ast}(A(Y), B), \quad \ast= 0, 1.
	\end{equation}
\end{defn}

This is a straightforward generalization of the notion \textit{topological $K$-group of $B$}: by definition, $K^\top_{\ast}(G; B) \coloneqq \KKH_{\ast}(G; \C, B)$ (\cite[Definition 9.1]{BCH}, \cite[page 9]{BMP}).

Any {\em proper} {\em $G$-compact} space $Y$ gives rise to a canonical element 
	\begin{equation}
	\lambda_{Y \rtimes G} \in K_{0}(\Co(Y) \rtimes G) = \KK_{0}(\C, \Co(Y) \r G)
	\end{equation}
by \cite[page 178]{KS2}.
\begin{defn}\label{defn assembly} Let $A$ be a \ca~and let $B$ be a \ga. The map
\begin{align}\label{assembly}
\beta^{A, B}_{G}:\KKH_{\ast}(G; A, B) \to \KK_{\ast}(A, B \r_\red G),
\end{align}
induced at the direct limit level by the composition
\begin{align*}
\beta^Y_G: \KK^G_{\ast}&(A(Y), B) \overset{j^G_\red}{\longrightarrow} \KK_{\ast}(A(Y) \r_\red G, B \r_\red G)\\
&= \KK_{\ast}(A \otimes (\Co(Y) \rtimes_{\red} G), B \rtimes_{\red} G) \overset{(1_{A} \x \lambda_{Y \rtimes G}) \x}{\longrightarrow} \KK_{\ast}(A, B \r_\red G),
\end{align*}
is called the (reduced) {\em naive assembly map} for $(A, B)$. Here $j^{G}_{\red}$ denote the reduced {\em descent map} of Kasparov (cf. \cite[3.11]{Ka}).
\end{defn}

\begin{conj}[The naive Baum-Connes Conjecture in $\KK$-theory]\label{naive BC for KK} Let $A$ be a \ca~ and let $B$ be a \ga. We say that $G$ satisfies the naive Baum-Connes conjecture for $(A, B)$ if the naive assembly map $\beta_{G}^{A, B}$ is an {\em isomorphism} of abelian groups. 
\end{conj}

The reason for the ``naiveness'' is that while the right-hand-side of the conjecture is $\sigma$-additive in the first variable, the left-hand-side is not. See Subsections \ref{sub nontrivial colimits} and \ref{sub open subgroups} for more details.

\begin{rem} \begin{enumerate}[(i)] 
\item The original conjecture of Baum and Connes states that for any group the assembly map is an isomorphism for the pair $(\C, \C)$. 
\item As stated in the introduction, counterexamples to the conjecture for $(\C, B)$ were constructed by Higson, Lafforgue and Skandalis \cite{HLS}.
\item Let 
	\begin{equation}\label{dimension-drop}
	\mathbb I_n \coloneqq \left\{f \in C([0,1], \mathbb M_n) \mid f(0) = 0, f(1) \in \C I\right\}, \quad n \in \Z_{\ge 1},
	\end{equation}
denote the $n$-th {\em dimension-drop algebra}. Then the mod-$n$ $K$-theory can be computed (cf.\ \cite{DL}) by 
	\begin{equation}
	K_{*}(D; \Z/n\Z) \cong \KK_{*}(\mathbb I_n, D).
	\end{equation}
Thus the Baum-Connes conjecture for $(\mathbb I_{n}, B)$ can be considered as a Baum-Connes conjecture for $B$ in mod-$n$ $K$-theory. 
\end{enumerate}
\end{rem}

\begin{rem}[Nontrivial action on $A$]
Suppose that $A$ is a \ga~ with a not necessarily trivial action of $G$. Then the topological $\KK$-groups of $(A, B)$ can be defined exactly as in Definition \ref{def top KK} and the definition of the assembly map can be modified to give an assembly map $\KKH_{*}(G; A, B) \to \KK_{*}(A^{G}, B \rtimes_{\red}G)$, where $A^{G}$ denote the fixed-point algebra of Kasparov \cite[Definition 3.2]{Ka}. However, the right-hand-side ``forgets'' too much information about the action of $G$ on $A$ for the assembly map to be an isomorphism in general. 

For instance, suppose that $G$ is a finite group. Then $\E{G} =\pt$ and $\KKH_{*}(G; A, B) = \KK^{G}_{*}(A, B)$ for any $(A, B)$. Let $H$ be an subgroup of $G$ and let $G$ act on $G/H$ by left-translation. Then 
	\begin{align}
	\KK^{G}_{*}(C(G/H), \C) \cong \KK^{H}_{*}(\C, \C) \cong \KK_{\ast}(\C, C^{*}(H)),
	\end{align}
by \cite[Proposition 5.14]{CE-P}, whereas
	\begin{align}
	\KK_{*}(C(G/H)^{G}, \C \r G) \cong \KK_{*}(\C, C^{*}(G)),
	\end{align}
since $C(G/H)^{G} \cong C(G\backslash G/H) \cong \C$. These can be quite different.
\end{rem}

\subsection{Compact groups}
 Let $G$ be a {\em compact} group. Then $\E{G} = \{\mathrm{pt}\}$ and $\lambda_{\{\mathrm{pt}\} \rtimes G} \in K_{0}(C^{*}(G))$ is the class of the central projection in $C^{*}(G)$  corresponding to the {\em trivial} representation of $G$. Let $A$ be a \ca~ (with the trivial $G$-action) and let $B$ be a \ga. The topological $\KK$-groups of $(A, B)$ are simply the equivariant $\KK$-groups:
	\[\KKH_{\ast}(G; A, B) = \KK^{G}_{\ast}(A, B)\]
and the assembly map equals the Green-Julg isomorphism 
	\begin{align}\label{Green-Julg map}
	\beta_{G}^{{A, B}}: KK^{G}_{\ast}(A, B) \overset{\cong}{\longrightarrow} \KK_{\ast}(A, B \rtimes G).
	\end{align}
See \cite[Proposition 6.25]{Tu-H} for more details. Hence we have the following.
\begin{prop}[Green-Julg Isomorphism]\label{compact case} Compact groups satisfy the naive Baum-Connes conjecture for any pair $(A, B)$. \qed
\end{prop}

\subsection{$\sigma$-additivity}\label{sub nontrivial colimits}
In this subsection, we explain why the Conjecture~\ref{naive BC for KK} is called naive. We claim that we have a ``problem'', whenever we have a ``nontrivial'' colimit in the definition of the naive topological $\KK$-group (cf.\ Definition \ref{defn assembly}). 
 
Indeed, let $A_{i}$ be \ca s, $i \ge 1$. Then
	 \begin{align}
	\KK_{\ast}(\oplus_{i}A_{i}, B \r_{\red} G) \cong \prod_{i} \KK_{\ast}(A_{i}, B \r_{\red} G),
	\end{align}
by the $\sigma$-additivity of $\KK$ in the first variable \cite[Theorem 2.9]{Ka}.
 On the other hand,
	\begin{align}
	\KKH_{\ast}(G; \oplus_{i}A_{i}, B) &= \colim_{\substack{Y \subseteq \E{G} \\ \text{$G$-compact}}} \KK^{G}_{\ast}((\oplus_{i}A_{i})(Y), B)\\
	&\cong \colim_{\substack{Y \subseteq \E{G} \\ \text{$G$-compact}}} \KK^{G}_{\ast}(\oplus_{i}A_{i}(Y), B)\\
	&\cong \colim_{\substack{Y \subseteq \E{G} \\ \text{$G$-compact}}} \prod_{i}\KK^{G}_{\ast}(A_{i}(Y), B),
	\end{align}
 again using \cite[Theorem 2.9]{Ka}. But this is {\em not} necessarily isomorphic to
 	\begin{align}
 	\prod_{i} \KKH_{\ast}(G; A_{i}, B) =  \prod_{i} \colim_{\substack{Y \subseteq \E{G} \\ \text{$G$-compact}}}\KK^{G}_{\ast}(A_{i}(Y), B), 
	\end{align}
since limits and colimits do {\em not} commute in general. Hence we cannot expect $G$ to satisfy the naive Baum-Connes conjecture for $(\oplus_{i}A_{i}, B)$ even if it does for $(A_{i}, B)$ for all $i \ge 1$.

\subsection{Ascending union of open  subgroups}\label{sub open subgroups}
We give an explicit example illustrating the difficulties of \ref{sub nontrivial colimits}.

Let $A$ be a \ca~ with the trivial action of $G$ and let $B$ be a \ga.
\begin{prop}[{cf.\ \cite[Theorem 5.1]{BMP}}]\label{BMP51}
Let $H$ be an {\em open} subgroup of $G$. Then the inclusion of $H$ in $G$ determines a homorphism of abelian groups 
	\[\KKH_{\ast}(H; A, B) \to \KKH_{\ast}(G; A, B).\]
If $G = \bigcup G_{n}$ is the union of ascending sequence of open subgroups 
	\[G_{1} \subseteq G_{2} \subseteq \dots \subseteq G,\]  
then
	\begin{equation}
	\KKH_{\ast}(G; A, B) \cong \colim_{n \to \infty} \KKH_{\ast}(G_{n}; A, B).
	\end{equation}	
\end{prop}
\begin{proof} The proof of \cite[Theorem 5.1]{BMP} applies ad verbatim, once we notice that since the $G$-action on $A$ is trivial, for any $H$-space $X$, we have 
	\[\Ind^{G}_{H}A(X) \cong A \otimes \Ind^{G}_{H} \Co(X),\]
where $\Ind^{G}_{H}: \KK^{H} \to \KK^{G}$ is the induction functor of Kasparov (cf.\ \cite[Theorem 3.5]{Ka}). 
\end{proof}

On the analytical side, we have the following. 
\begin{prop}[{cf.\ \cite[Theorem 4.1]{BMP}}]\label{BMP41} Let $H$ be an open subgroup of $G$, then canonical inclusion $C_{c}(H, B) \to C_{c}(G, B)$ extends to an injective $\ast$-homomorphism
	\[B \r_{\red} H \to B \r_{\red} G.\]
If $G = \bigcup G_{n}$ is the union of ascending sequence of open subgroups, then	
	\begin{equation}
	B \r_{\red} G \cong \colim_{n \to \infty} B \r_{\red}G_{n}.
	\end{equation}
\end{prop}
\begin{proof} See the proof of \cite[Theorem 4.1]{BMP}. 
\end{proof}

\begin{defn}\label{defn KK-compact} We say that a \ca\ $A$ is $\KK$-compact if $\KK_{\ast}(A, -)$ is \cont, {\em i.e.}\ commutes with colimits.
\end{defn}

\begin{ex} 
\begin{enumerate}
\item If $A$ satisfies the UCT (cf.\ Theorem~\ref{UCT}) and has finitely generated $K$-theory, then $A$ is $\KK$-compact (cf.\ \cite[Theorem 7.13]{RS}).  In particular, the dimension-drop algebra $\mathbb I_{n}$ of (\ref{dimension-drop}) is $\KK$-compact.
\item If $A$ has a $K$-amenable Poincar{\'e} dual in the sense of \cite[VI.4.$\beta$]{Co}, then $A$ is $\KK$-compact.
\end{enumerate}
\end{ex}

\begin{thm}[{cf.\ \cite[Theorem 6.3]{BMP}}] Let $A$ be a $\KK$-compact \ca~ and let $B$ be a \ga. Suppose that $G$ is the union of ascending sequence of open subgroups $G_{n}$, each satisfying the naive Baum-Connes conjecture (\ref{naive BC for KK}) for $(A, B)$. Then $G$ satisfies the naive Baum-Connes conjecture for $(A, B)$. 
\end{thm}
\begin{proof} Since $\KK_{\ast}(A, -)$ is continuous, 
	\[\KK_{\ast}(A, B \r_{\red} G) \cong \colim_{n \to \infty} \KK_{\ast}(A, B \r_{\red} G_{n}).\]
	Now the proof of \cite[Theorem 6.3]{BMP} applies.
\end{proof}

Since $\KK$-theory is {\em not} continuous in the second variable (cf.\ \cite[19.7.2]{Ba}), we {\em cannot} expect $\KK_{\ast}(A, B\r_{\red}G)$ to be isomorphic to $\colim_{n \to \infty} \KK_{\ast}(A, B \r_{\red}G_{n})$ without restrictions on $A$. We demonstrate by example that the continuity of $\KK_{\ast}(A, -)$ is necessary. This particular example was suggested by Nigel Higson (in the context of subsection \ref{sub nontrivial colimits}).

\begin{ex}\label{naive BC not good ex}
Let $G$ denote the (discrete) abelian group $\bigoplus_{k \ge 1} \Z/2\Z$ and let $G_{n} \coloneqq \bigoplus_{k = 1}^{n} \Z/2\Z$ considered as a subgroup of $G$. Then $G = \bigcup_{n \ge 1} G_{n}$. Note that abelian groups satisfy the Baum-Connes conjecture for any $(\C, B)$.

Let $A \coloneqq c_{0}(\Lambda)$ for some countable set $\Lambda$ and let $B \coloneqq \C$. Then $B \r G_{n} = C^{*}(G_{n}) \cong C^{*}(\Z/2\Z)^{\otimes n} \cong (\C^{2})^{\otimes n}$ and the inclusion map $B \r G_{n} \to B \r G_{n+1}$ is given by $f \mapsto f \otimes 1_{\C^{2}}$. Hence $\KK(\C, C^{*}(G_{n})) \cong (\Z^{2})^{\x n} \cong \Z^{2^{n}}$ and the map induced by the inclusion is given by $\Z^{2^{n}} \ni p \mapsto (p, p) \in \Z^{2^{n+1}}$.

On the topological side, by Proposition \ref{BMP51},
	\begin{align*}
	\KKH(G; c_{0}(\Lambda), \C) &= \colim_{n \to \infty} \KKH(G_{n}; c_{0}(\Lambda), \C)\\
	&= \colim_{n \to \infty} \KK(c_{0}(\Lambda), \C \r G_{n})\\
	&= \colim_{n \to \infty} \prod_{\Lambda} \KK(\C, C^{*}(G_{n}))\\
	&= \colim_{n \to \infty} \prod_{\Lambda} \Z^{2^{n}}. 
	\end{align*}
On the analytical side, by Proposition \ref{BMP41}, 
	\begin{align*}
	\KK(c_{0}(\Lambda), \C \r_{\red} G) &= \prod_{\Lambda} \KK(\C, \C \r G)\\
	&= \prod_{\Lambda} \colim_{n \to \infty} \KK(\C, \C \r G_{n})\\
	&= \prod_{\Lambda} \colim_{n \to \infty} \KK(\C, C^{*}(G_{n}))\\
	&= \prod_{\Lambda} \colim_{n \to \infty} \Z^{2^{n}}.
	\end{align*}
Now it is a simple algebraic exercise to show that the two groups are different. Hence $G = \bigoplus_{k \ge 1} \Z/2\Z$ does {\em not} satisfy the naive Baum-Connes conjecture for $(A, B) = (c_{0}(\Lambda), \C)$.
\end{ex}

\subsection{Continuity}\label{subsec Continuity}  Now we show that if $A$ is $\KK$-compact, then the naive Baum-Connes conjecture is stable under taking inductive limits of $G$-algebras. See Corollary \ref{BC direct sum} for the precise statement.

\begin{lem}[{cf.\ \cite[Subsection 1.1]{CEO2}}]\label{KK is GD} Let $A$ and $B$ be \ga s. Then
	\begin{equation}
	\F_H(\Co(Y)) \coloneqq \KK^{H}(\Res^H_G A(Y), \Res^H_GB), \quad H \in \S(G)
	\end{equation}
is a Going-Down functor in the sense of \cite[Definitions 1.1]{CEO2}.

If the $G$-action on $A$ is trivial, then 
	\begin{equation}
	\F^{n}(G) = \KKH_{n}(G; A, B).	
	\end{equation}
\end{lem}
\begin{proof}
The functor $\F^\ast_H$ is homotopy invariant by \cite[Proposition 2.5]{Ka} and satisfies the Restriction axiom by \cite[Theorem 5.8]{Ka}. Let 
	\[0 \rightarrow \Co(U) \rightarrow \Co(Y) \rightarrow \Co(Y \backslash U) \rightarrow 0\]
be a short exact sequence of proper commutative $H$-algebras. Then it is equivariantly semi-split by \cite[Corollary 6.2]{KS1} and hence so is the sequence 
	\[0 \rightarrow A(U) \rightarrow A(Y) \rightarrow A(Y \backslash U) \rightarrow 0.\]
Then the same corollary  implies that $\F^\ast_H$ is half-exact. Thus $\F$ satisfies the Cohomology axioms. Finally, by \cite[Lemma 3.6]{Ka} there is a natural isomorphism: 
	\[A(\Ind^G_H(\Co(Y))) \cong \Ind^G_H(A(Y))\]
and \cite[Remark 5.4]{KS1} (see \cite[Proposition 5.14]{CE-P} for a slightly more general version) proves that $\F$ satisfies the Induction axiom. The last statement is clear.
\end{proof}

Often we omit the restriction functors from notation.

\begin{lem}[{cf.\ \cite[Proposition 7.1]{CE-P}}]\label{cont of KKH}
Let $A$ be a $\KK$-compact algebra. Then the functor $\KKH_{\ast}(G; A, -)$ is continuous. 
\end{lem}
\begin{proof} Let $B = \colim_{i \to \infty} B_i$ be a direct limit of \ga s 
	\[ \dots \to B_i \to B_{i+1} \to \dots.\]

For $H \in \S(G)$, set
	\begin{align*}
	\F^{\ast}_H(\Co(Y)) &= \colim_{i \to \infty} \KK^{H}_{\ast}(A(Y), B_i) ~{\rm and}\\ 
	\G^{\ast}_H(\Co(Y))&= \KK^{H}_{\ast}(A(Y), B).
	\end{align*} 
Then $\F$ and $\G$ are \GD functors in the sense of \cite[Definition 1.1]{CEO2} and the natural maps $\Lambda_H: \F_H \rightarrow \G_H$, induced by $B_i \rightarrow B$, form a \GD transformation in the sense of \cite[Definition 1.3]{CEO2}. Moreover, in the notation of \cite[Section 1]{CEO2}, 
	\begin{align*}
	\F^{n}(G) &= \colim_{\substack{Y \subseteq \E{G} \\ \text{$G$-compact}} } \colim_{i \to \infty} \KK^{G}_{n}(A(Y), B_i)\\
	 &\cong \colim_{i \to \infty}\colim_{\substack{Y \subseteq \E{G} \\ \text{$G$-compact}} }  \KK^{G}_{n}(A(Y), B_i)\\
	 & \cong \colim_{i \to \infty} \KKH_{n}(G; A, B_{i})\qquad\text{and}\\
	\G^{n}(G) & = \KKH_{n}(G; A, B).
	\end{align*}
We need to show that $\Lambda^{n}(G): \F^{n}(G) \rightarrow \G^{n}(G)$ is an isomorphism.
	
Let $V$ be a finite dimensional Euclidean space equipped with a linear action of $K$. We have natural isomorphisms 
	\begin{equation}
	\KK^K(A(V), B_{(i)}) \cong \KK^K(A, B_{(i)}(V)) \cong \KK(A, B_{i}(V) \r K)
	\end{equation}
by Kasparov's Bott periodicity theorem (cf.\ \cite[Lemma 7.7~]{CE-T}) and the Green-Julg theorem (cf.\ Proposition \ref{compact case}). Since $B(V) \r K \cong \colim_{i \to \infty} (B_i(V) \r K$, the following commutative diagram 	\begin{equation}
	\xymatrix{
	\colim_{i \to \infty} \KK^K_{\ast}(A(V), B_i) \ar[r]^{\qquad\Lambda_{K}} \ar[d]^{\cong} & \KK^K_{\ast}(A(V), B) \ar[d]^{\cong}\\
	\colim_{i \to \infty} \KK_{\ast}(A, B_i(V) \r K)  \ar[r]^{\qquad\cong} & \KK_{\ast}(A, B(V) \r K)}
	\end{equation}
proves that the map $\Lambda_K: \F_K(\Co(V)) \to \G_K(\Co(V))$ is an isomorphism. Thus by \cite[Theorem 1.4]{CEO2}, $\Lambda^{{n}}(G)$ is an isomorphism and $\KKH(G; A, -)$ is continuous.
\end{proof}

\begin{cor}[{cf.\ \cite[Proposition 2.5]{CEN}}]\label{BC direct sum} Let $A$ be a $\KK$-compact algebra and let $ \dots \to B_{i} \to B_{i+1} \to \dots$ be an inductive system of \ga s. Suppose that either $G$ is exact {\em or} all the connecting maps $B_{i} \to B_{i+1}$ are injective. Then if $G$ satisfies the naive Baum-Connes conjecture for $(A, B_{i})$ for all $i$, then it satisfies for $(A, \colim_{i} B_{i})$. \qed
\end{cor}

\section{The Baum-Connes Conjecture for $\KK$-theory, attempt $2$}
In this section, we consider an alternative generalization of the Baum-Connes conjecture to $\KK$-theory. This generalization is already considered in \cite{Ka}, in the case of almost connected groups.

\subsection{Almost connected groups}
A topological group is said to be {\em almost connected} if its group of connected components is compact. The Baum-Connes conjecture for almost connected groups is known for the pair $(\C, \mathcal{K})$ with any action on $\mathcal{K}$, where $\mathcal{K}$ is the algebra of compact operators on a separable Hilbert space (cf. \cite{CEN}). 

In this section, we study the naive Baum-Connes conjecture for an almost connected group $G$ for a general pair $(A, B)$. This will serve as a toy model and leads to the second approach to the Baum-Connes conjecture for $\KK$-theory. The following two characteristics make it particularly nice to work with:
\begin{enumerate}[(a)]
\item it admits a $G$-compact universal proper space (hence difficulties from subsection~\ref{sub nontrivial colimits} do not arise)
\item it has a $\gamma$-element (cf. \cite[Definition 1.7]{CE-P}).
\end{enumerate}

Let $K$ be a maximal compact subgroup of $G$. Then the quotient $X \coloneqq G/K$ equipped with the left-translation action of $G$ is a universal proper $G$-space by \cite[Main Theorem]{Ab}.  Since $X$ is $G$-compact, we have 
	\begin{equation}
	\KKH_{\ast}(G; A, B) = \KK^{G}_{\ast}(A(X), B).
	\end{equation}
Let $\P  = C_{\tau}(X)$ denote the graded algebra of the $C_{0}$-sections of the Clifford bundle on $X$ and let $d = [d_{X}] \in \KK^{G}_{0}(\P, \C)$ denote the Dirac element of $X$ (cf.\ \cite[Definition-Lemma 4.2]{Ka}).

\begin{thm}\label{H = MN}
Let $G$ be an almost connected group and let $A$ be a \ca~ and let $B$ be a \ga. 
Then the assembly map $\beta^{A, B}_{G}$ can be identified with the ``multiplication by the Dirac element''  
	\[\x j^{G}_{\red}(1_{B} \x d): \KK_{\ast}(A, (B \x \P) \r_{\red} G) \to \KK_{\ast}(A, B \r_{\red} G)\]
via the commutative diagram
	\begin{equation}\label{assembly = multiplication}\xymatrix{
	\KKH_{\ast}(G; A, B \x \P) \ar[r]^{\beta^{A, B \x \P}_{G}}_{\cong} \ar[d]^{ \x (1_{B} \x d)}_{\cong} &\KK_{\ast}(A, (B \x \P) \r_{\red} G) \ar[d]^{\x j^{G}_{\red}(1_{B} \x d)}\\
	\KKH_{\ast}(G; A, B) \ar[r]^{\beta^{A, B}_{G}} &\KK_{\ast}(A, B \r_{\red} G). 
	}\end{equation}
\end{thm}
This is certainly well-known to the experts, but since we could not find any direct reference, we provide a proof (compare \cite[Theorem 5.10]{Ka}). First we fix some notation.

\begin{notation} Let $G$ be a topological group and let $H$ be a closed subgroup. For an $H$-algebra $D$, the canonical Morita equivalence from $D \r_{\red} H$ to $(\Ind_{H}^{G}D) \r_{\red}G$ is denoted by $x_{D}$ (cf.\ \cite[Theorem 3.15]{Ka}). For a $G$-algebra $E$, the canonical $G$-isomorphism from $E(G/H)$ to $\Ind^{G}_{H}\Res^{H}_{G} E$, given by $\Co(G/H, E) \ni f \mapsto [\widetilde f: g \mapsto g f(g^{-1}K)] \in \Ind^{G}_{H}\Res^{H}_{G} E$ is denoted by $\varphi_{E}$ (cf.\ \cite[Lemma 3.6]{Ka}).
\end{notation}

\begin{lem}[{cf.\ \cite[Proposition 2.3]{CE-P}}]\label{lem shapiro for acg} Let $G$ be an almost connected group and let $K \subseteq G$ be a maximal compact subgroup and let $X = G/K$. Let $A$ be a \ca~ and let $D$ be a $K$-algebra. Then the following diagram is commutative:
	\begin{equation}\label{shapiro for acg}\xymatrix{
	\KK^{K}_{\ast}(A, D) \ar[r]^{\beta^{A, D}_{K}}_{} \ar[d]^{\Ind^{G}_{K}}_{} &\KK_{\ast}(A, D \r K) \ar[dd]^{\cdot \x [x_{D}]}_{}\\
	\KK^{G}_{\ast}(\Ind^{G}_{K}A, \Ind^{G}_{K} D) \ar[d]^{[\varphi_{A}] \x \cdot}_{}  & \\
	\KK^{G}_{\ast}(A(X), \Ind^{G}_{K} D) \ar[r]^{\beta^{A, \Ind^{G}_{K} D}_{G}} &\KK_{\ast}(A, (\Ind^{G}_{K} D) \r_{\red} G). 
	}\end{equation}
\end{lem}
\begin{proof}
Take $x \in \KK^{K}_{\ast}(A, D)$. Then we need to show that
	\begin{equation}\label{shapiro for acg, eqn}(1_{A} \x \lambda_{X \r G}) \x j^{G}_{\red}([\varphi_{A}] \x \Ind^{G}_{K} x) = (1_{A} \x \lambda_{\pt \r K}) \x j^{K}_{\red} x \x [x_{D}].
	\end{equation}

Since the action on $A$ is trivial, $A \r K \cong A \otimes \C \r K$ and $\Ind^{G}_{K}A \cong A \otimes \Ind^{G}_{K} \C$ and under this identification $x_{A} = 1_{A} \x x_{\C}$ and $\varphi_{A} = 1_{A} \x \varphi_{\C}$. Thus, (\ref{shapiro for acg, eqn}) is the consequence of the following identities:
\begin{enumerate}
\item 	$\lambda_{X \rtimes G} \x j^{G}_{\red}[\varphi_{\C}] = \lambda_{\{\mathrm{pt}\} \rtimes K} \x [x_{\C}]$ (cf.\ \cite[(2.4)]{CE-P}) and
\item $j^{K}_{\red}x \x [x_{D}] = [x_{A}] \x j^{G}_{\red}\Ind^{G}_{K}x.$ (cf.\ \cite[Corollary 3.15]{Ka}).
\end{enumerate}	
\end{proof}

\begin{lem}\label{lem induction} Let $G$ be an almost connected group and let $K \subseteq G$ be a maximal compact subgroup. Let $A$ be a \ca~ and let $E$ be a \ga. Then the induction map
	\[\KK^{K}_{\ast}(A, \Res^{K}_{G} E) \to \KK^{G}_{\ast}(\Ind^{G}_{K}A, \Ind^{G}_{K}\Res^{K}_{G}E)\]
is an isomorphism.
\end{lem}
See Corollary~\ref{cor induction plus} for a stronger statement.
\begin{proof}
We keep the notation $X = G/K$. Consider the diagram, 
	\begin{equation}
	 \xymatrix{
	\KK^{K}_{*}(A, \Res^{K}_{G}E) \ar[rr]^{\cdot \x 1_{\Co(X)}} \ar[dd]^{\Ind_{K}^{G}}& &\KK^{K}_{*}(A(X), \Res^{K}_{G}E(X))\\
	& \KK^{G}_{*}(A, E) \ar[ul]^{\Res^{K}_{G}} \ar[rd]^{\cdot \x 1_{\Co(X)}}&\\
	\KK^{G}_{*}(\Ind_{K}^{G}A, \Ind_{K}^{G}\Res^{K}_{G}E) \ar[rr]^{[\varphi_{A}] \otimes \cdot \x [\varphi_{E}^{-1}]} & & \KK^{G}_{*}(A(X), E(X)) \ar[uu]^{\Res^{K}_{G}}.
	}\end{equation}

The lower-left corner is commutative by \cite[Theorem 3.6]{Ka} and the upper-right corner is commutative by the functoriality of the   restriction map. Moreover, the restriction map $\Res^{K}_{G}: \KK^{G}_{*}(A, E) \to \KK^{K}_{*}(A, \Res^{K}_{G}E)$ is surjective by \cite[Corollary 5.7]{Ka}, therefore the rectangle on the outside is commutative. 

The classes $[\varphi_{A}]$ and $[\varphi_{E}]$ are equivalences by construction and the restriction $\Res^{K}_{G}: \KK^{G}_{*}(A(X), E(X)) \to \KK^{K}_{*}(A(X), \Res^{K}_{G}E(X))$ is an isomorphism by \cite[Theorem 5.8]{Ka}. It remains to show that $\cdot \x 1_{\Co(X)}: \KK^{K}_{*}(A, D) \to \KK^{K}_{*}(A(X), D(X))$ is an isomorphism. This follows from the equivariant Bott periodicity of Kasparov (cf.\ \cite[Lemma 7.7]{CE-T}) since, by \cite[Corollary A.6]{Ab}, $X$ can be given the structure of a real vector space for which the action of $K$ is linear. 
\end{proof}

\begin{proof}[Proof of Theorem~\ref{H = MN}]
Commutativity of (\ref{assembly = multiplication}) follows from the multiplicative property of $j^{G}_{\red}$ (cf.\ \cite[Theorem 3.11]{Ka}). The vertical map on the left 
	\[\x (1_{B} \x d): \KK^{G}_{\ast}(A(X), B \x \P) \to \KK^{G}_{\ast}(A(X), B)\]
is an isomorphism by \cite[Theorem 5.8]{Ka}, with inverse $\cdot \x (1_{B} \x \eta)$, where $\eta = \eta_{X} \in \KK^{G}_{0}(\C, \P)$ is the {\em dual-Dirac} element of Kasparov (cf.\ \cite[Definition-Lemma 5.1]{Ka}). Finally, let $C_{V}$ denote the Clifford algebra of the cotangent space to $X = G/K$ at $K \in X$ (cf.\ \cite[Theorem 5.10]{Ka}). Then we have
	\begin{equation}\label{P is compactly induced}
	\P  = \Ind^{G}_{K}\Res^{K}_{G} (C_{V}).
	\end{equation}
Combining Lemmas \ref{lem shapiro for acg} and \ref{lem induction} with the Green-Julg isomorphism (Proposition \ref{compact case}), we see that 	the assembly map $\beta^{A, B \x \P}_{G}$ is an isomorphism. This completes the proof.
\end{proof}

\begin{cor}\label{cor GJ for ACG} Let $G$ be an almost connected group and let $A$ be a \ca~ and let $B$ be a \ga. Then the assembly map gives an isomorphism 
 	\begin{equation}
	\KKH(G; A, B) \cong \KK(A, B \r_{\red} G) \x j^{G}_{\red}(1_{B} \x \gamma),
	\end{equation}
where $\gamma = \gamma_{G} \coloneqq \eta \x d \in \KK^{G}(\C, \C)$ is the $\gamma$-element of Kasparov (cf.\ \cite[Theorem 5.7]{Ka}).

In particular, $G$ satisfies the Baum-Connes conjecture for $(A, B)$ if and only if $j^{G}_{\red}(1_{B} \x \gamma)$ acts as the identity on $\KK(A, B \r_{\red} G)$.
\end{cor}

The right-hand-side of the expression is called the {\em $\gamma$-part} of $\KK(A, B \r_\red G)$.
\begin{proof} This is a well-known argument. It follows from the proof of Theorem~\ref{H = MN} that for any $x \in \KK^{G}_{\ast}(A(X), B)$
	\begin{align*}
	\beta^{A, B}_{G}(x) &= \beta^{A, B \x \P}(x \x (1_{B} \x \eta)) \x j^{G}_{\red}(1_{B} \x d)\\
	& = (1_{A} \x \lambda_{X \r G}) \x j^{G}_{\red}(x \x (1_{B} \x \eta)) \x j^{G}_{\red}(1_{B} \x d)\\
	& = (1_{A} \x \lambda_{X \r G}) \x j^{G}_{\red}(x)  \x j^{G}_{\red}(1_{B} \x \eta d)\\
	& = \beta^{A, B}_{G}(x) \x j^{G}_{\red}(1_{B} \x \gamma). 
	\end{align*}
The proof is completed using the identity $\gamma^{2} = \gamma$.	
\end{proof}

\begin{rem} This corollary is not necessarily true for other groups with a $\gamma$-element in the sense of \cite[Definition 1.7]{CE-P}, see Example~\ref{naive BC not good ex}.
\end{rem}

\subsection{Strong Baum-Connes conjecture for almost connected groups}

Applying Yoneda's lemma to Theorem~\ref{H = MN}, we get the following.
\begin{cor}\label{cor strong BC} Let $G$ be an almost connected group and let $B$ be a \ga. Then $G$ satisfies the Baum-Connes conjecture for $(A, B)$ for all \ca s $A$ if and only if $j^{G}_{\red}(1_{B} \x d) \in \KK((B \x P) \r_{\red} G, B \r_{\red} G)$ is invertible if and only if $j^{G}_{\red}(1_{B} \x \gamma) = 1_{B \r_{\red} G} \in \KK(B \r_{\red} G, B \r_{\red} G)$.
\qed
\end{cor}

If $G$ and $B$ satisfies the equivalent properties of Corollary~\ref{cor strong BC}, we say that $G$ satisfies the {\em strong} Baum-Connes conjecture for $B$ (cf.\ \cite[Definition 9.1]{MN-T}).

\begin{ex}
Any almost connected group with the Haagerup property satisfies $\gamma = 1 \in \KK^{G}(\C, \C)$ (cf.\ \cite{HK}), thus satisfies the strong Baum-Connes conjecture for any \ga~ $B$. Examples include $\mathrm{SO}(n, 1)$ and $\mathrm{SU}(n, 1)$.
\end{ex}

\begin{cor}[{cf.\ \cite[Proposition 9.5]{MN-T}}]\label{B type I} Let $G$ be an almost connected group and let $B$ be a type $I$ \ga. Then $G$ satisfies the strong Baum-Connes conjecture for $B$ if and only if $G$ satisfies the Baum-Connes conjecture for $(\C, B)$ and $B \r_{\red} G$ satisfies the UCT.
\end{cor}
We include the short proof for the convenience of the reader.
\begin{proof} 
First note that the algebra $(B \x P)\r_{\red} G$ satisfies the UCT, since it is Morita equivalent to $\Res^{K}_{G}(B \x C_{V}) \r K$, which is type $I$ by Takesaki's theorem (\cite[Theorem 6.1]{Ta}). Now suppose that $G$ satisfies the strong Baum-Connes conjecture. Then, clearly, $G$ satisfies the Baum-Connes conjecture for $(\C, B)$ and $B \r_{\red} G$ satisfies the UCT by virtue of being $\KK$-equivalent to $(B \x P)\r_{\red} G$. 

Conversely, suppose that $G$ satisfies the Baum-Connes conjecture for $(\C, B)$ and $B \r_{\red} G$ satisfies the UCT. Then by \cite[Proposition 23.10.1]{Ba},
	\[j^{G}_{\red}(1_{B} \x d) \in \KK((B \x P) \r_{\red} G, B \r_{\red} G)\]
is invertible.
\end{proof}

\begin{exs}\label{strong BC examples}
\begin{enumerate}
\item\label{strong BC for almost connected} Any almost connected group satisfies the strong Baum-Connes conjecture for $B = \mathcal{K}$. Indeed, let $G$ be almost connected. Then $G$ satisfies the Baum-Connes conjecture for $(\C, \mathcal{K})$ by \cite{CEN} and $\mathcal{K} \r_{\red} G$ satisfies the UCT by \cite[Proposition 5.1]{CEO2}.   
\item\label{strong BC for Sp} Let $\Gamma$ be a discrete subgroup of $\mathrm{Sp}(n, 1)$ of finite covolume, $n \ge 2$. Then $B = \Ind^{G}_{\Gamma} \C$ is commutative (hence type $I$) but $(\Ind^{G}_{\Gamma}\C) \r_{\red} G$ does not satisfy the UCT. Indeed by \cite[Corollaire 4.2]{Sk88}, the algebra $C^{\ast}_{\red}\Gamma$, which is Morita equivalent to $(\Ind^{G}_{\Gamma}\C) \r_{\red} G$, is not even $\KK$-equivalent to a nuclear algebra, let alone an abelian one. Hence $\mathrm{Sp}(n, 1)$  do {\em not} satisfy the strong Baum-Connes conjecture for $\Ind^{G}_{\Gamma}\C$.  On the other hand, $\mathrm{Sp}(n, 1)$ does satisfy the usual Baum-Connes conjecture for $(\C, B)$ for any $B$ (cf.\ \cite{Ju02}). This example is due to Skandalis \cite{Sk88}.
\end{enumerate}  
\end{exs}

It follows from the equation (6.1) of \cite{CE-P}, 
	\begin{equation}
	j^{G}_{\red}(1_{\Ind^{G}_{K}D} \x \gamma_{G}) = [x_{D}^{-1}] \x j^{K}_{\red}(1_{D} \x \gamma_{K}) \x [x_{D}] = 1_{\Ind^{G}_{K}D \r_{\red}G},
	\end{equation}
that $G$ satisfies the strong Baum-Connes conjecture for the induced algebra $\Ind^{G}_{K}D$, for any $K$-algebra $D$ (See also \cite[Proposition 10.1]{MN-T}).
This allows us improve on Lemma~\ref{lem induction}.
\begin{cor}\label{cor induction plus} Let $G$ be an almost connected group and let $K \subseteq G$ be a maximal compact subgroup. Let $A$ be a \ca~ and let $D$ be a $K$-algebra. Then the induction map
	\[\Ind^{G}_{K}: \KK^{K}_{\ast}(A, D) \to \KK^{G}_{\ast}(\Ind^{G}_{K}A, \Ind^{G}_{K}D)\]
is an isomorphism.
\end{cor}
\begin{proof} Every map except $\Ind^{G}_{K}$ in the commutative diagram (\ref{shapiro for acg}) is an isomorphism, hence so is $\Ind^{G}_{K}$.
\end{proof}

\subsection{The Baum-Connes conjecture for $\KK$-theory} 
When $G$ is almost connected, Theorem~\ref{H = MN} can be used to prove many nice properties of the assembly map. For the general case, we turn around everything, and reformulate the conjecture so that Theorem~\ref{H = MN} becomes a tautology. 

We recall some terminology from \cite{MN-T}. From now on, we work with equivariance with respect to transformation groupoids. This generality is needed in Section~\ref{sec Comparison}: we deduce Corollary~\ref{main thm for CI}, which is used in the proof of the Comparison Theorem~\ref{comparison},  from the forgetful isomorphism of Theorem~\ref{thm forgetful}.

Let $X$ be a \gs.

\begin{defn}[{\cite[Definition 4.1]{MN-T}}] An \xga~ is {\em compactly induced} if it is isomorphic to $\Ind^{G}_{K} D$ for some compact subgroup of $K \subseteq G$ and some $K$-algebra $D$. Let $\CI \subseteq \KK^{X \r G}$ denote the full subcategory of compactly induced algebras and let $\LCI$ denote the localizing subcategory generated by $\CI$. A morphism $f \in \KK^{X \r G}$ is called a {\em weak equivalence} if $\Res^{K}_{G}f \in \KK^{X \r K}$ is an isomorphism for all compact subgroups $K \subseteq G$.
\end{defn}

\begin{defn}[{\cite[Definition 4.5]{MN-T}}]\label{defn Dirac el} An element $d \in \KK^{X \r G}(\P, \Co(X))$ is called a {\em Dirac morphism} for $X \r G$ if $d$ is a {\em $\CI$-simplicial approximation} of $\C \in \KK^{G}$, that is,  
	\begin{enumerate}
	\item $\P$ is an object of $\LCI$ and
	\item  $d$ is a weak equivalence.
	\end{enumerate}
\end{defn}

By \cite[Proposition 4.6]{MN-T}, Dirac morphisms exist, uniquely up to isomorphism, for any transformation groupoid. It follows that $\LCI$ is a {\em coreflective} subcategory of $\KK^{\XG}$.

\begin{ex} Let $G$ be an almost connected group and let $K$ be a maximal compact subgroup. Then the Dirac element $d = d_{G/K} \in \KK^{G}(\P, \C)$ of \cite[Definition-Lemma 4.2]{Ka} is a Dirac morphism for $G$ in the sense of Definition~\ref{defn Dirac el} (Strictly speaking we need to replace $P$ by an ungraded algebra.) 
\end{ex}

Let $d \in \KK^{X \r G}(\P, \Co(X))$ be a Dirac morphism for $X \r G$. 
\begin{defn} Let $A$ be a \ca.\ and let $B$ be an \xga. We define the {\em topological $\KK$-group} of $(A, B)$ as
	\begin{equation}
	\KKMN_{\ast}(X \r G; A, B) \coloneqq \KK_{\ast}(A, (B \x_{X} P) \r_{\red} G)
	\end{equation}
and the (reduced) {\em assembly map} as
	\begin{equation}
	\mu^{A, B}_{X \r G} \coloneqq \cdot \x j^{G}_{\red}(1_{B} \x_{X} d): \KKMN_{\ast}(X \r G; A, B) \to \KK_{\ast}(A, B \r_{\red} G).
	\end{equation}
\end{defn}
Theorem 5.2 of \cite{MN-T} shows that this is indeed a generalization of the Baum-Connes conjecture. We write $K^{\top}_{\ast}(X \r G; B)$ for $\KKMN_{\ast}(X \r G; \C, B)$.

\begin{rem}\label{rem shapiro KKMN} By \cite[Lemma 5.1]{MN-T}, if $d \in \KK^{G}(P, \C)$ is a Dirac morphism for $G$ then $p_{X}^{*}(d) \in \KK^{X \r G}(P(X), \Co(X))$,  where $p_{X}: X \to \pt$,  is a Dirac morphism for $X \r G$ and the natural identification 
	\begin{equation}
	\mathcal F: B \x_{X} P(X) \cong B \x P
	\end{equation}
satisfies $\mu_{X \r G} = \mu_{G} \circ \mathcal F_{\ast}$. 
\end{rem}

\begin{conj} (The Baum-Connes Conjecture in KK-theory).\label{BCMN} Let $A$ be a \ca\ and let B be a $G$-algebra. We say that $G$ satisfies the Baum-Connes conjecture for $(A, B)$ if the assembly map $\mu^{A, B}_{G}$ is an isomorphism of abelian groups. 
\end{conj}

This formulation doesn't have the shortcoming of the naive version, described in Subsection~\ref{sub nontrivial colimits}: If $G$ satisfies the Baum-Connes conjecture for $(A_{i}, B)$ for all $i$, then it satisfies for $(\oplus_{i} A_{i}, B)$.

\section{Comparison of the two approaches}\label{sec Comparison}

We know that the two formulations of the generalized Baum-Connes conjecture are not equivalent.
\begin{ex} Let $G \coloneqq \bigoplus_{k \ge 1} \Z/2\Z$ and $A = c_{0}(\Z)$ and $B = \C$. Then $G$ satisfies Conjecture~\ref{BCMN} for $(A, B)$ by the $\sigma$-additivity of $\KK$ in the first variable, but not Conjecture~\ref{naive BC for KK} as demonstrated in Example~\ref{naive BC not good ex}.
\end{ex} 

\subsection{The comparison map}
First we generalize the naive topological $\KK$-theory to transformation groupoids, following \cite{Tu-H} and \cite{CEO1}. Let $X$ be a $G$-space. 

\begin{defn} Let $A$ be a \ca\ and let $B$ be an \xga. We define the naive topological $\KK$-groups as 
	\begin{equation}
	\KKH_{\ast}(X \r G; A, B) \coloneqq \colim_{\substack{Y \subseteq X \t \E{G} \\ \text{$G$-compact}}} \KK^{X \r G}_{\ast}(A(Y), B).
	\end{equation}
\end{defn} 
As in \cite[Section 1]{CEO1}, there is a {\em forgetful} map 
	\begin{equation}\label{forgetful map}
	\mathcal F: \KKH_{\ast}(X \r G; A, B) \to \KKH_{\ast}(G; A, B)
	\end{equation} and an assembly map 
	\begin{equation}
	\beta^{A, B}_{X \r G}: \KKH_{\ast}(X \r G; A, B) \to \KK_{\ast}(A, B \r_{\red}G)
	\end{equation}
satisfying $\beta_{X \r G} = \beta_{G} \circ \mathcal F$, defined inductively via maps
	\begin{align}
	\mathcal F_{Y}&: \KK^{X \r G}(A(Y), B) \overset{F_{X}}{\to} \KK^{G}(A(Y), B) \overset{(\pi_{2}|_{Y})_{\ast}}\to \KK^{G}(A(\pi_{2}(Y)), B)\\
	\beta^{Y}_{X \r G}&: \KK^{X \r G}(A(Y), B) \overset{F_{X}}{\to} \KK^{G}(A(Y), B) \overset{\beta^{Y}_{G}}{\to} \KK(A, B \r_{\red} G),
	\end{align}
where $F_{X}: \KK^{X \r G} \to \KK^{G}$ is the forgetful map and $\pi_{2}: X \t \E{G} \to \E{G}$ is the projection onto the second coordinate. 

Now we define a comparison map
	\begin{equation}
	\nu^{A, B}_{X \r G}: \KKH_{\ast}(X \r G; A, B) \to \KKMN_{\ast}(X \r G; A, B).	\end{equation}
Let $d \in \KK^{G}(\P, \Co(X))$ be a Dirac morphism for $X \r G$. Then we have a commutative diagram 
	\begin{equation}\label{comparison}\xymatrix{
	\KKH_{\ast}(X \r G; A, B \x_{X} \P) \ar[r]^{\beta^{A, B \x_{X} \P}_{X \r G}} \ar[d]^{ \cdot \x (1_{B} \x_{X} d)} & \KKMN_{\ast}(X \r G; A, B) \ar[d]^{\mu^{A, B}_{X \r G}}\\
	\KKH_{\ast}(X \r G; A, B) \ar[r]^{\beta^{A, B}_{X \r G}} &\KK_{\ast}(A, B \r_{\red} G). 
	}\end{equation}

\begin{lem}\label{lem mult by w.e.} Let $x \in \KK^{X \r G}(B, B')$ be a weak equivalence and let $A$ be a \ca. Then the natural map 
	\[\cdot \x x: \KKH_{\ast}(X \r G; A, B) \to \KKH_{\ast}(X \r G; A, B')\]
an isomorphism.
\end{lem}
\begin{proof} Let $Y \subseteq X \t \E{G}$ be a $G$-compact subset. Then $Y$ is proper and, by \cite[Corollary 7.3]{MN-T},  the algebra $A(Y)$ belongs to $\LCI$ (using $G$-compactness, we see that $A(Y)$ is contained in the triangulated subcategory generated by $\CI$). By \cite[Proposition 4.4]{MN-T}, 
	\begin{equation}
	\cdot \x x: \KK^{X \r G}(A(Y), B) \to \KK^{X \r G}(A(Y), B')	
	\end{equation}
is an isomorphism.
\end{proof}

As a corollary, the leftmost vertical map
	\begin{equation}
	\cdot \x (1_{B} \x_{X}d):\KKH_{\ast}(X \r G; A, B \x_{X} P) \to \KKH_{\ast}(X \r G; A, B)
	\end{equation}
is an isomorphism.

\begin{defn}
We define the {\em comparison} map as the composition
	\begin{equation}
	\nu^{A, B}_{X \r G} \coloneqq \beta^{A, B \x_{X} \P}_{X \r G} \circ (\cdot \x (1_B \x_X d))^{-1},
	\end{equation}
going from $\KKH_{\ast}(X \r G; A, B)$ to $\KKMN_{\ast}(X \r G; A, B)$.
\end{defn}	
This is an analogue of the map $\nu$ of \cite[Section 5]{Tu-H}. It follows from the commutativity of (\ref{comparison}) that 
	\begin{equation}
	\mu^{A, B}_{X \r G} \circ \nu^{A, B}_{X \r G} = \beta^{A, B}_{X \r G}.
	\end{equation}

Our main theorem is the following.
\begin{thm}[Comparison]\label{thm comparison} Let $A$ be a $\KK$-compact algebra (cf.\ Definition~\ref{defn KK-compact}) and $B$ be an \xga. Then the comparison map $\nu^{A, B}_{X \r G}$ is an isomorphism.
\end{thm}
The main difficulty in the proof is that we do not know if we can choose $P$, the source of the Dirac morphism, to be a proper algebra. However, this is the case for $G$ almost connected and this fact turns out to be sufficient.

\subsection{Proof of the Comparison Theorem \ref{thm comparison}}
First we suppose that $X \r G$ has a Dirac morphism $d \in \KK^{X \r G}(\P, \Co(X))$ with $\P$ proper, that is, $P$ admits a $\Co(X \t \E{G})$-structure. 

For any $G$-invariant subsets $V \subseteq Y \subseteq X \t \E{G}$ with $V$ open and $Y$ $G$-compact, we have the descent isomorphism of Kasparov and Skandalis 
	\begin{equation}\label{descent map}
	\KK^{Y \r G}_{\ast}(A(Y), B \x_{X} P_{V}) \cong \KK_{\ast}(A, (B \x_{X} P_{V}) \r G),
	\end{equation}
which is given by the forgetful map $F_{Y}: \KK^{Y \r G} \to \KK^{X \r G}$ followed by the assembly map $\beta^{Y}_{X \r G}$ (cf.\ \cite[Proposition 6.25]{Tu-H}). Here $P_{V} = \Co(V)P$ is the restriction to $V$. 

Let $i_{V}: P_{V} \to P$ denote the inclusion and let $d_{V} = [i_{V}] \x d$. Then we have a natural map
	\begin{align}
	\KK_{*}(A, (B \x_{X} P_{V}) &\r G) \cong \KK^{Y \r G}_{*}(A(Y), B \x_{X} P_{V}) \overset{F_{Y}}{\to} \KK^{X \r G}_{*}(A(Y), B \x_{X} P_{V})\\
	&\overset{\cdot \x d_{V}}{\to} \KK^{X \r G}_{*}(A(Y), B) \to \KKH_{*}(X \r G; A, B).
	\end{align}

If $A$ is $\KK$-compact, then taking the colimit over $V$ (cf.\ \cite[Proposition 5.7]{Tu-H}), we get a map
	\begin{equation}
	\kappa^{A, B}_{X \r G}: \KKMN_{*}(X \r G; A, B) = \KK_{*}(A, (B \x_{X}P) \r G) \to \KKH_{*}(X \r G; A, B).
	\end{equation}	

It is clear that 
	\begin{align}
	\beta^{A, B}_{X \r G} \circ \kappa^{A, B}_{X \r G} &= \mu^{A, B}_{X \r G}\quad\text{and}\\ 
	\nu^{A, B}_{X \r G} \circ \kappa^{A, B}_{X \r G} &= \Id.
	\end{align} 

\begin{prop} Suppose that $X \r G$ has a Dirac morphism $d \in \KK^{X \r G}(\P, \Co(X))$ with $\P$ proper. Let $A$ be a $\KK$-compact algebra and let $B$ be an \xga. Then the comparison map $\nu^{A, B}_{X \r G}$ is an isomorphism, with inverse $\kappa^{A, B}_{X \r G}$.
\end{prop} 
\begin{proof} 
We need to show that $\kappa^{A, B}_{X \r G} \circ \nu^{A, B}_{X \r G} = \Id$. By \cite[Corollary 7.2]{MN-T}, the element
	\begin{equation}p^{*}_{\E{G}}(d) \in \KK^{(X \t \E{G})\r G}(\P(\E{G}), \Co(X \t \E{G}))
	\end{equation}
 is invertible, where $p_{\E{G}}: X \t \E{G}   \to X$ is projection onto the first coordinate. Let 
 	\begin{equation}\theta \in \KK^{(X \t \E{G})\r G}(\P(\E{G}), \Co(X \t \E{G}))	\end{equation}
denote the inverse. First we claim that the inverse of $\cdot \x (1_{B} \x_{X} d)$ is given by multiplication by $\theta$ on the left. More explicitly, let $j_{Y}: Y \subseteq X \t \E{G}$ be the inclusion of a $G$-compact subset and let $\theta_{Y} \coloneqq j_{Y}^{*}(\theta) \in \KK^{Y \r G}(\Co(Y), P \x_{X} \Co(Y))$. Let
	\begin{align}
	 F_{Y}&: \KK^{Y \r G} \to \KK^{X \r G} 
	\end{align}
denote the forgetful map. 

Let $x \in \KK^{X \r G}(A(Y), B)$. Then
	\begin{align*}
	(\theta \x x) \x d &= (1_{A} \x F_{Y}\theta_Y \x_{A(Y) \x_{X} P} (x \x_{X} 1_{P})) \x_{B \x_{X} P} (1_{B} \x_{X} d)\\
	&= F_{Y}\theta_Y \x_{\Co(Y)} (x \x_{X} d) \qquad\text{(cf.\ \cite[Lemme 5.5]{Tu-H})}\\
	&=F_{Y}\theta_Y \x_{\Co(Y)} (d \x_{X} x)\\
	&=F_{Y}\theta_Y \x_{\Co(Y)} (d \x_{X} 1_{\Co(Y)}) \x x\\
	&=F_{Y}j_{Y}^{*}(\theta \x p_{\E{G}}^{*}d ) \x x\\
	&= x.
	\end{align*}
Now let $x \in \KK^{X \r G}(A(Y), B)$. Then $\nu^{A, B}_{X \r G}(x) = \beta^{A, B \x_{X} P}_{X \r G}(\theta \x x)$ and we need to write it in a form pluggable to $\kappa^{A, B}_{X \r G}$. 
 
The descent isomorphism and the continuity of $K$-theory imply that 
	\begin{equation}
	\KK^{Y \r G}(\Co(Y), P \x_{X} \Co(Y)) \cong \colim_{V}\KK^{Y \r G}(\Co(Y), P_{V} \x_{X} \Co(Y)).
	\end{equation}
Consequently, there exists $V \subseteq X \t \E{G}$ open and $\theta_{Y, V} \in \KK^{Y \r G}(\Co(Y), P_{V} \x_{X} \Co(Y))$ such that 
	\begin{equation}
	\theta_{Y} = \theta_{Y, V} \x_{P_{V}} [i_{V}].
	\end{equation}
Moreover, according to \cite[Proposition 5.12]{Tu-H}, there exists a $G$-compact subset $L \subseteq X \t \E{G}$, containing both $V$ and $Y$, and $\theta' \in \KK^{L \r G}(\Co(L), P_{V} \x_{X} \Co(Y))$, where $\Co(L)$ acts on the first component of $P_{V} \x_{X} \Co(Y)$, such that 
	\begin{equation}
	F_{L}\theta' = [j_{Y, L}] \x F_{Y}\theta_{Y, V}
	\end{equation}
in $\KK^{X \r G}(\Co(L), P_{V} \x_{X} \Co(Y))$, where $j_{Y, L}: Y \to L$ is the inclusion.

Then, as in \cite{Tu-H}, we can write 
	\begin{align*}
	\beta^{A, B \x_{X} P}_{X \r G} (\theta \x x) & = \beta^{Y} (F_{Y}\theta_{Y} \x_{\Co(Y)} x)\\
	& = \beta^{Y}((F_{Y}\theta_{Y, V} \x_{P_{V}} [i_{V}])\x_{\Co(Y)} x)\\
	& = \beta^{Y}(F_{Y}\theta_{Y, V} \x_{\Co(Y) \x_{X} P_{V}} ([i_{V}]\x_{X} x))\\
	& = \beta^{Y}(F_{Y}\theta_{Y, V} \x_{\Co(Y) \x_{X} P_{V}} (x \x_{X} [i_{V}]))\\
	& = \beta^{Y} (F_{Y}\theta_{Y, V} \x_{\Co(Y)} x) \x j^{G}[i_{V}]\\
	& = \beta^{L} (F_{L}(\theta' \x_{\Co(Y)} x)) \x j^{G}[i_{V}].
	\end{align*}	
It follows that 
	\begin{align*}
	\kappa^{A, B}_{X \r G} (\nu^{A, B}_{X \r G}(x)) & = \kappa^{A, B}_{X \r G}(\beta^{A, B \x_{X} P}_{X \r G} (\theta \x x))\\
	& = F_{L}(\theta' \x_{\Co(Y)} x) \x_{P_{V}} d_{V}\\
	& = [j_{Y, L}] \x F_{Y}\theta_{Y, V} \x (x \x_{X} d_{V})\\
	& = [j_{Y, L}] \x F_{Y} \theta_{Y, V} \x [i_{V}] \x d \x x\\
	& = [j_{Y, L}] \x (F_{Y}\theta_{Y} \x_{P} d) \x x\\
	& = [j_{Y, L}] \x x. 
	\end{align*}
This completes the proof.
\end{proof}

If $G$ is an almost connected group, then $X \r  G$ has a Dirac morphism with proper $P$. Thus, as in Corollary~\ref{cor GJ for ACG}, we get the following.
\begin{cor}\label{cor forgetful ACG} Let $G$ be an almost connected group. Let $A$ be a $\KK$-compact \ca~ and let $B$ be an \xga. Then 
	\[\beta^{A, B}_{X \r G}: \KKH_{*}(X \r G; A, B) \to \KK_{*}(A, B\r_{\red} G)\]
is an isomorphism onto the $\gamma$-part of $\KK_{*}(A, B\r_{\red} G)$.
\qed
\end{cor}

As a corollary, we get the following.
\begin{thm}[{Forgetful Isomorphism, cf.\ \cite[Theorem 0.1]{CEO1}}]\label{thm forgetful} Let $A$ be a $\KK$-compact algebra. Then the forgetful map
	\[\mathcal F: \KKH_{*}(X \r G; A, B) \to \KKH_{*}(G; A, B)\]
of (\ref{forgetful map}) is an isomorphism.
\end{thm}
\begin{proof} Corollary~\ref{cor forgetful ACG} implies that the theorem holds for $G$ almost connected. Now the proof of \cite[Theorem 0.1]{CEO1} applies.
\end{proof}
 
\begin{cor}\label{main thm for CI} Let $A$ be a $\KK$-compact algebra. Then the assembly map $\beta^{A, B}_{G}$ is an isomorphism for $B \in \CI$.
\end{cor}
\begin{proof} Proceeds as in \cite[Section 4]{CEO1}.
\end{proof}

Now we are ready to prove the Comparison Theorem~\ref{thm comparison}.
\begin{proof}[Proof of the Comparison Theorem \ref{thm comparison}]
We need to show that $\nu^{A, B}_{X \r G}$, or equivalently $\beta^{A, B \x_{X} P}_{X \r G}$, is an isomorphism. By Theorem~\ref{thm forgetful}, it is enough to consider the case $X = \pt$. Let  $\BC(G; A)$ denote the full subcategory of $G$-algebras $E \in \KK^{G}$ such that $\beta^{A, E}_{G}$ is an isomorphism. Then $\BC(G; A)$ is clearly a triangulated subcategory of $\KK^{G}$. Moreover, since $A$ is $\KK$-compact, $\BC(G; A)$ is closed under countable direct sums by Corollary~\ref{BC direct sum} and contains $\CI$ by Corollary~\ref{main thm for CI}. Hence $\LCI \subseteq \BC(G; A)$. Now it is enough to notice that $B \x P$ belongs to $\LCI$ (cf.\ \cite[Lemma 4.2]{MN-T}).
\end{proof}

\section{The Universal Coefficient Theorem}\label{section UCT}

In this section we develop a Universal Coefficient Theorem (UCT) for topological $\KK$-functors and prove Theorem \ref{mod n}. As an application, we get an alternative proof of Theorem~\ref{thm comparison}, in the case $A$ satisfies the UCT and has finitely generated $K$-theory (such $A$'s are $\KK$-compact). 

First we recall the UCT of Rosenberg and Schochet (\cite[Section IX.23]{Ba}).
\begin{thm}[UCT \cite{RS}]\label{UCT} A \ca~ $A$ is $\KK$-equivalent to an abelian \ca~  if and only if it satisfies the UCT for every $B$, that is, there is a natural short exact sequence:
	\begin{equation}
	\xymatrix{
	\Ext^{*}_\Z(K_\ast(A), K_\ast(B)) \ar@{>->}[r] &\KK_\ast(A, B) \ar@{->>}[r] & \Hom^{*}_\Z(K_\ast(A), K_\ast(B)).
	}\end{equation}
\qed	
\end{thm}
In this situation, we simply say that {\em $A$ satisfies the UCT}. The full subcategory of $\KK$ of algebras satisfying the UCT is the localizing subcategory $\langle \C \rangle \subset \KK$ generated by $\C$ (cf.\ \cite[Section 2.5]{MN-T}).

As in \cite{CEO2}, we develop an abstract UCT  first and specialize it to the topological $\KK$-functors. 

\begin{defn} Let $\mathcal C \subseteq \KK$ be a triangulated subcategory containing $\C$. A {\em UCT functor} on $\mathcal C$ is a cohomological functor $F: \mathcal C \rightarrow \Ab$, to the category of abelian groups, equipped with a zero-graded natural transformation 
	\begin{equation}
	\gamma_{A}: F_\ast(A) \rightarrow \Hom^{\ast}_\Z(K_\ast(A), F_\ast(\C)),
	\end{equation}
such that $\gamma_{A}$ is an isomorphism whenever $K_\ast(A)$ is free and finitely generated. If, in addition, $\mathcal C$ is localizing and $\gamma_{A}$ is an isomorphism whenever $K_{\ast}(A)$ is free, then we say that $F$ is {\em $\sigma$-UCT}.
\end{defn}

\begin{prop}[Abstract UCT]\label{abstract UCT} Let $\mathcal C \subseteq \KK$ be a triangulated subcategory containing $\C$, and let $F$ be a UCT functor on $\mathcal C$. Then for every \ca\ $A$ in $\mathcal C$ with finitely generated $K$-theory, there is a natural short exact sequence, called the UCT exact sequence:
	\begin{equation}
	\xymatrix{
	\Ext^{*}_\Z(K_\ast(A), F_\ast(\C)) \ar@{>->}[r] &F_\ast(A) \ar@{->>}[r] & \Hom^{*}_\Z(K_\ast(A), F_\ast(\C)).
	}\end{equation}
If $F$ is $\sigma$-UCT, then the UCT exact sequence exists for all $A$ in $\mathcal C$ (with no restriction on $K_{\ast}(A)$).
\end{prop}
This is standard, but we include a proof here, because the proof of the usual UCT in \cite{Ba} uses an injective resolution of $K_{\ast}(B)$, whereas we use a free resolution of $K_{\ast}(A)$. As usual, it is enough to assume that $F$ is defined only on $\ast$-homomorphisms, not arbitrary $\KK$-morphisms.
\begin{proof} We proceed as in \cite[Section 3]{CEO2}. In both cases, it follows from Schochet's construction of the geometric resolution (cf.\ \cite[Proposition 23.5.1]{Ba}) that there exists an algebra $R$ in $\mathcal C$ and a $*$-homomorphism $\varphi: R \to A \x \mathcal K$, where $\mathcal K$ is the algebra of compact operators on a separable Hilbert space, such that $K_{\ast}(R)$ is free and $\varphi_{\ast}: K_{\ast}(R) \to K_{\ast}(A \x \mathcal K) \cong K_{\ast}(A)$ is surjective. The rotated mapping cone triangle
	\begin{equation}
	\Sigma R \overset{-\Sigma\varphi}{\to} \Sigma(A \x \mathcal K) \to C_{\varphi} \overset{\varepsilon}{\to} R
	\end{equation}
is an exact triangle in $\mathcal C$. This gives a free resolution
	\begin{equation}
	0 \to K_{\ast}(C_{\varphi}) \overset{K_{\ast}(\varepsilon)}{\to} K_{\ast}(R) \to K_{\ast}(A) \to 0
	\end{equation}
of $K_{\ast}(A)$ and consequently 
	\begin{align}
	\Hom_{\Z}(K_{\ast}(A), F_{\ast}(\C)) &\cong \ker \Hom_{\Z}(K_{\ast}(\varepsilon), F_{\ast}(\C))\qquad \text{and}\\
	\Ext_{\Z}(K_{\ast}(A), F_{\ast}(\C)) &\cong \coker \Hom_{\Z}(K_{\ast}(\varepsilon), F_{\ast}(\C)).
	\end{align}
Moreover, since we have a commutative diagram
	\begin{equation}
	\xymatrix{
	F_{\ast}(R) \ar[r]^{F_{\ast}(\varepsilon)} \ar[d]^{\gamma_{R}}_{\cong} & F_{\ast}(C_{\varphi}) \ar[d]^{\gamma_{C_{\varphi}}}_{\cong}\\
	\Hom_{\Z}(K_{\ast}(R), F_{\ast}(\C)) \ar[r]&	\Hom_{\Z}(K_{\ast}(C_{\varphi}), F_{\ast}(\C)),
	}
	\end{equation}
we may identify
	\begin{align}
	\label{identify gamma1}\Hom_{\Z}(K_{\ast}(A), F_{\ast}(\C)) &\cong \ker F_{\ast}(\varepsilon)\qquad \text{and}\\
	\label{identify gamma2}\Ext_{\Z}(K_{\ast}(A), F_{\ast}(\C)) &\cong \coker F_{\ast}(\varepsilon).
	\end{align}
Finally, since $F$ is a cohomological functor, we have a short exact sequence
	\begin{equation}
	0 \to \coker F_{\ast}(\varepsilon) \to F_{\ast}(\Sigma(A \x \mathcal K)) \to \ker F_{\ast}(\Sigma\varepsilon) \to 0,
	\end{equation}
which in combination with the identifications (\ref{identify gamma1}) and (\ref{identify gamma2}) completes the proof.
\end{proof}

For a fixed \ca\ $B$, the functor $A \mapsto \KK(A, B)$ is a $\sigma$-UCT functor on $\langle \C \rangle$. Applying the Abstract UCT we get Theorem~\ref{UCT}. As a corollary, we obtain the following. 
\begin{thm}\label{UCT for KKtop} Let $B$ be a \ga. For any algebra $A$ satisfying the UCT, we have the following natural short exact sequences and the assembly maps induce a map of short exact sequences
	\begin{equation*}
	\xymatrix{
	\Ext^{*}_\Z(K_\ast(A), K^{\top}_\ast(G; B)) \ar@{>->}[r] \ar[d] &\KKMN_\ast(G; A, B) \ar@{->>}[r] \ar[d] & \Hom^{*}_\Z(K_\ast(A), K_\ast^{\top}(G; B)) \ar[d]\\
	\Ext^{*}_\Z(K_\ast(A), K_{\ast}(B \r_{\red} G)) \ar@{>->}[r] &\KK_\ast(A, B \r_{\red} G) \ar@{->>}[r] & \Hom^{*}_\Z(K_\ast(A), K_\ast(B \r_{\red} G)).
	}\end{equation*}
\end{thm}
\begin{proof}  Follows from the functoriality of the UCT sequence: the assembly map $\cdot \x j^{G}_{\red}(1_{B} \x d)$ induces a map of short exact sequences between the UCT sequences for $(A, (B \x P) \r_{\red} G)$ and $(A, B \r_{\red} G)$. 
\end{proof}

Applying the Five-Lemma, we obtain the following.
\begin{cor}\label{UCT and BCMN} Let $B$ be a \ga. Suppose that $G$ satisfies the Baum-Connes conjecture for $(\C, B)$. Then for any algebra $A$ satisfying the UCT, $G$ satisfies the Baum-Connes conjecture for $(A, B)$. 
\qed
\end{cor}

Next we consider the UCT for naive $\KK$-theory.
\begin{thm}\label{UCT for naive KK} Let $B$ be a \ga. Then for any $A$ satisfying the UCT and having finitely generated $K$-theory, we have the following natural short exact sequences and the assembly maps induce a map of short exact sequences
	\begin{equation*}
	\xymatrix{
	\Ext^{*}_\Z(K_\ast(A), K^{\top}_\ast(G; B)) \ar@{>->}[r] \ar[d] &\KKH_\ast(G; A, B) \ar@{->>}[r] \ar[d] & \Hom^{*}_\Z(K_\ast(A), K_\ast^{\top}(G; B)) \ar[d]\\
	\Ext^{*}_\Z(K_\ast(A), K_{\ast}(B \r_{\red} G)) \ar@{>->}[r] &\KK_\ast(A, B \r_{\red} G) \ar@{->>}[r] & \Hom^{*}_\Z(K_\ast(A), K_\ast(B \r_{\red} G)).
	}\end{equation*}
\end{thm}
\begin{proof}
Let $\mathcal C$ denote the full subcategory of $\langle \C \rangle$ consisting of algebras with finitely generated $K$-theory.  It is clear that $\mathcal C$ is a triangulated subcategory containing $\C$. Let $B$ be a \ga. We consider the functor $F: \mathcal C \to \Ab$ given by
	\begin{equation}
	F(A) \coloneqq \KKH(G; A, B)
	\end{equation}
on objects. If $x $ is a morphism in $\mathcal C(A', A)$, then it can be considered an element of $\KK^{G}_{\ast}(A', A)$ naturally and $F(x): F(A') \to F(A)$ is given by the multiplication 
	\begin{equation}
	x \x \cdot: \KK^{G}(A(Y), B) \to \KK^{G}(A'(Y), B)
	\end{equation}
at the inductive limit level. Then $F$ is a cohomological functor on $\mathcal C$. 
Moreover, using the identity $K_{\ast}(A) = \KK_{\ast}(\C, A)$, we get a map
	\begin{equation}
	\gamma_{A}: \KKH_{\ast}(G; A, B) \to \Hom_{\Z}^{\ast}(K_{\ast}(A), \KKH_{\ast}(G; \C, B)).
	\end{equation}	
This is certainly a natural transformation and we need to show that if $K_{\ast}(A)$ is finitely generated and free then $\gamma_{A}$ is an isomorphism. Using the finite-additivity of both sides, it is enough to consider the cases $A = \C$ and $A = \Sigma \C$, which are obvious.

The last assertion is clear.	
\end{proof}

We note that since $\KKH(G; A, B)$ is not necessarily $\sigma$-additive in $A$, the functor $F$ above is {\em not} $\sigma$-UCT in general.

Applying the Five-Lemma, we get the following.
\begin{cor}\label{UCT and naive BC} Let $B$ be a \ga. Suppose that $G$ satisfies the Baum-Connes conjecture for $B$. Then for any $A$ satisfying the UCT and having finitely generated $K$-theory, $G$ satisfies the naive Baum-Connes conjecture (\ref{naive BC for KK}) for $(A, B)$. \qed
\end{cor}

Let $A$ be an algebra satisfying the UCT and having finitely generated $K$-theory. Then the comparison map $\nu^{A, B}_{G}$ is an isomorphism. Indeed, by Corollary~\ref{UCT and naive BC}, it is enough to show that $G$ satisfies the usual Baum-Connes conjecture for $B \x P$. But this is clear since $B \x P$ belongs to $\LCI$ by \cite[Lemma 4.2]{MN-T} and elements of $\LCI$ satisfy the usual Baum-Connes conjecture by \cite[Theorem 5.2]{MN-T}.

In particular, the two versions of the mod-$n$ Baum-Connes conjecture are equivalent and they are implied by the usual Baum-Connes conjecture.

\bibliographystyle{amsalpha}
\bibliography{master}

\end{document}